 \documentclass[english]{bourbaki}

\newcommand{\SL}{\operatorname{SL}}
\newcommand{\BD}{\operatorname{BD}}

\newcommand{\D}{\operatorname{D}}
\newcommand{\Ext}{\operatorname{Ext}}
\newcommand{\Fix}{\operatorname{Fix}}
\newcommand{\Jac}{\operatorname{Jac}}
\newcommand{\NonSing}{\operatorname{NonSing}}
\newcommand{\Spec}{\operatorname{Spec}}
\newcommand{\Stab}{\operatorname{Stab}}
\newcommand{\Supp}{\operatorname{Supp}}
\newcommand{\age}{\operatorname{age}}
\newcommand{\card}{\operatorname{card}}
\newcommand{\codim}{\operatorname{codim}}
\newcommand{\diag}{\operatorname{diag}}
\renewcommand{\div}{\operatorname{div}}
\newcommand{\wt}{\operatorname{wt}}
\newcommand{\Hilb}{\operatorname{Hilb}}
\newcommand{\AHilb}{A\hbox{-}\!\Hilb}
\newcommand{\GHilb}{G\hbox{-}\!\Hilb}
\newcommand{\SU}{\operatorname{SU}}
\newcommand{\aff}{\mathbb A}
\newcommand{\LL}{\mathbb L}
\newcommand{\PP}{\mathbb P}
\newcommand{\C}{\mathbb C}
\newcommand{\Q}{\mathbb Q}
\newcommand{\R}{\mathbb R}
\newcommand{\Z}{\mathbb Z}
\newcommand{\sD}{\mathcal D}
\newcommand{\sF}{\mathcal F}
\newcommand{\sI}{\mathcal I}
\newcommand{\Oh}{\mathcal O}
\newcommand{\sV}{\mathcal V}
\newcommand{\sZ}{\mathcal Z}
\newcommand{\bm}{\mathbf m}
\newcommand{\dd}{\mathrm d}
\newcommand{\stringy}{_{\mathrm{string}}}
\newcommand{\open}{^\circ}

\newcommand{\iso}{\cong}
\newcommand{\into}{\hookrightarrow}
\newcommand{\broken}{\dasharrow}
\newcommand{\bij}{\leftrightarrow}
\newcommand{\wave}{\widetilde}
\newcommand{\blob}{{\scriptscriptstyle\bullet}}
\newcommand{\hidot}{^{\blob}}
\newcommand{\wD}{\wave D}

\newcommand{\Span}[1]{\left<#1\right>}
\newcommand{\dbrk}[1]{[\kern-.15em[#1]\kern-.15em]}
\newcommand{\al}{\alpha}
\newcommand{\be}{\beta}
\newcommand{\ga}{\gamma}
\newcommand{\ze}{\zeta}
\newcommand{\la}{\lambda}
\newcommand{\si}{\sigma}
\newcommand{\om}{\omega}
\newcommand{\De}{\Delta}
\newcommand{\Om}{\Omega}
\newcommand{\ep}{\varepsilon}
\newcommand{\fie}{\varphi}
\newcommand{\1}{^{-1}}
\newcommand{\tensor}{\otimes}
\numberwithin{equation}{section}

 \theoremstyle{plain}
 \newtheorem{princ}[defi]{Principle}
 \newtheorem{lem}[defi]{Lemma}
 \newtheorem{prodf}[defi]{Proposition--Definition}

 \makeatletter
 \newenvironment{pf}[1]{\par \normalfont
 \topsep6\p@\@plus6\p@ \trivlist \itemindent\normalparindent
 \item[\hskip\labelsep\normalfont\itshape 
 #1\@addpunct{\pointrait}]\ignorespaces}
 {\endtrivlist}
 \makeatother

\title{La correspondence de McKay}
\author{Miles Reid}
\curraddr{Math Inst.,\\ Univ. of Warwick, Coventry CV4 7AL}
\email{Miles@Maths.Warwick.Ac.UK} 
\date{}

 \bbknumero{867}
 \bbkannee{52\`eme ann\'ee, session de novembre 1999}

\begin{document}
\maketitle

\section*{Abstract:} Let $M$ be a quasiprojective algebraic manifold with $K_M=0$ and $G$
a finite auto\-morphism group of $M$ acting trivially on the canonical
class $K_M$; for example, a subgroup $G\subset\SL(n,\C)$ acting on $\C^n$
in the obvious way. We aim to study the quotient variety $X=M/G$ and its
resolutions $Y\to X$ (especially under the assumption that $Y$ has
$K_Y=0$) in terms of $G$-equivariant geometry of $M$. At present we know 4
or 5 quite different methods of doing this, taken from string theory,
algebraic geometry, motives, moduli, derived categories, etc.

For $G\subset\SL(n,\C)$ with $n=2$ or $3$, we obtain several methods of
cobbling together a basis of the homology of $Y$ consisting of algebraic
cycles in one-to-one correspondence with the conjugacy classes or the
irreducible representations of $G$.

 \section{Comment c'est}
 \label{sec:intro}

 \subsection{Model case: the binary dihedral group $\BD_{4n}$}
 \label{ssec:model} For $G\subset\SL(2,\C)$ a finite group, the quotient
variety $X=\C^2/G$ is called a {\em Klein quotient singularity}. I draw
the quotient map $\pi\colon\C^2\to X$ and the minimal resolution of
singularities $Y\to X$ together in the diagram:
 \begin{equation}
 \renewcommand{\arraystretch}{1.3}
 \begin{matrix} && \C^2 \\
 &&{\hphantom{\scriptstyle\pi}}\big\downarrow{\scriptstyle\pi} \\
 Y & \xrightarrow{\fie} & X
 \end{matrix}
 \notag
 \end{equation}
This situation has been well studied, since Klein around 1870 and Coxeter
and Du Val in the 1930s: the subgroup $G$ is classified as cyclic, binary
dihedral or a binary group corresponding to one of the Platonic solids; the
quotient singularity is a hypersurface $X\subset\C^3$ with defining
equation one of a list of simple functions. The resolution $Y$ is a
surface with $K_Y=\fie^*K_X$, and the exceptional locus $\fie\1(0)\subset
Y$ of the resolution consists of a bunch of $-2$-curves $E_i$ (that is,
$E_i\iso\PP^1_\C$ and $E_i$ has self-intersection $E_i^2=-2$), and the
intersection $E_iE_j$ is given by one of the Dynkin diagrams $A_n$, $D_n$,
$E_6$, $E_7$, $E_8$. To avoid writing out lists, let us simply discuss the
binary dihedral group
 \begin{equation}
 G=\BD_{4n}=\Span{\al,\be},\quad\text{where}\quad
 \al=\begin{pmatrix} \ep&0 \\ 0&\ep\1
 \end{pmatrix}, \quad
 \be=\begin{pmatrix} 0&1\\-1&0
 \end{pmatrix}
 \notag
 \end{equation}
where $\ep=\exp\frac{2\pi i}{2n}$. If $u,v$ are coordinates on $\C^2$, the
$G$-invariant polynomials are
 \begin{equation}
 \C[x,y,z]/(z^2-yx^2+4y^{n+1}), \quad\text{where}\quad x=u^{2n}+v^{2n},
 y=u^2v^2, z=uv(u^{2n}-v^{2n});
 \notag
 \end{equation}
thus the quotient variety is the singularity
$X:(z^2=yx^2-4y^{n+1})\subset\C^3$ of type $D_{n+2}$, and the quotient
morphism $\pi\colon(u,v)\mapsto(x,y,z)$. The resolution of singularities
$Y\to X$ has exceptional locus consisting of $-2$-curves
$E_1,\dots,E_{n+2}$ forming the $D_{n+2}$ configuration:
 
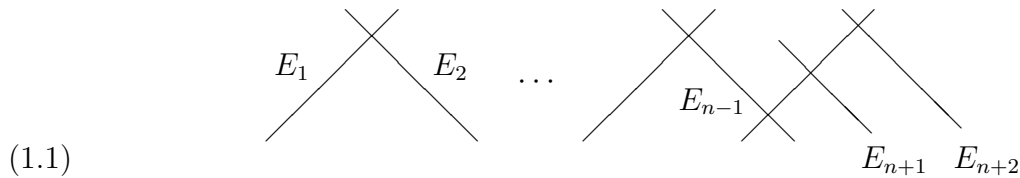
\begin{equation}
 \begin{picture}(80,70)
 \put(-87,35){$E_1$}
 \put(-90,10){\line(1,1){50}}
 \put(-27,35){$E_2$}
 \put(-60,60){\line(1,-1){50}}
 \put(5,30){$\cdots$}
 \put(30,10){\line(1,1){50}}
 \put(66,22){$E_{n-1}$}
 \put(60,60){\line(1,-1){50}}
 \put(90,10){\line(1,1){50}}
 \put(135,0){$E_{n+1}$}
 \put(104,48){\line(1,-1){35}}
 \put(170,0){$E_{n+2}$}
 \put(128,60){\line(1,-1){45}}
 \end{picture}
 \label{eq:resoln}
 \end{equation}

The classical McKay correspondence begins in the late 1970s with the
observation that the same graph arises in connection with the
representation theory of $G$. For a group $G$ and a given representation
$Q$, the {\em McKay graph} (or McKay quiver) has a node for each
irreducible representation, and an edge $V\to V'$ whenever $V'$ is a
direct summand of $V\tensor Q$. In our case, $\BD_{4n}$ has the
2-dimensional representations
 \begin{equation}
V_i\iso\C^2, \quad\text{with action}\quad
\al=\begin{pmatrix} \ep^i&0 \\ 0&\ep^{-i}
 \end{pmatrix}, \quad
\be=\begin{pmatrix} 0&1\\(-1)^i&0
 \end{pmatrix}\quad\text{for $i=0,\dots,n$.}
 \notag
 \end{equation}
This is irreducible for $0<i<n$, and splits into 2 eigenlines when $i=0$
or $n$. The inclusion $G\subset\SL(2,\C)$ provides the {\em given}
representation $Q=V_1$. It is a simple exercise \cite{Hwk} to write down
the action of $G$ on a basis $\{e_i\tensor e'_j\}$ of $Q\tensor V_i$ to
get $V_i\tensor Q=V_{i-1}\oplus V_{i+1}$ for $0<i<n$, so that the McKay
graph of $\BD_{4n}$ is the extended Dynkin diagram $\wD_{n+2}$:
 
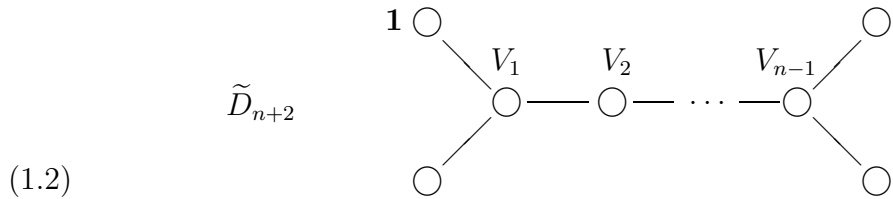
\begin{equation}
 \begin{picture}(150,75) 
 \put(-70,28){$\wD_{n+2}$}
 \put(-10,60){$\mathbf1$}
 \put(0,60){$\bigcirc$}
 \put(30,38){\line(-1,1){18}}
 \put(0,0){$\bigcirc$}
 \put(30,27){\line(-1,-1){18}}
 \put(30,30){$\bigcirc$}
 \put(44,33){\line(1,0){23}}
 \put(70,30){$\bigcirc$}
 \put(84,33){\line(1,0){15}}
 \put(101,32){ \dots }
 \put(124,33){\line(1,0){15}}
 \put(140,30){$\bigcirc$}
 \put(152,38){\line(1,1){18}}
 \put(170,60){$\bigcirc$}
 \put(152,27){\line(1,-1){18}}
 \put(170,0){$\bigcirc$}
 \put(30,45){$V_1$}
 \put(72,45){$V_2$}
 \put(130,45){$V_{n-1}$}
 \end{picture}
 \label{eq:repns}
 \end{equation}
Here $\mathbf1$ is the trivial 1-dimensional representation.

This example, and the other $\SL(2,\C)$ cases observed by McKay, suggest
that there is a one-to-one correspondence between the components of the
exceptional locus of $Y\to X$ in (\ref{eq:resoln}) and the nontrivial
irreducible representations of $G\subset\SL(2,\C)$ in (\ref{eq:repns}).
This talk explains this coincidence in several different ways, and
discusses higher dimensional generalisations.

 \subsection{General assumption}\label{ssec:ass}
I use the following diagram throughout:
 \begin{equation}
 \renewcommand{\arraystretch}{1.3}
 \begin{matrix} && M \\
 &&{\hphantom{\scriptstyle\pi}}\big\downarrow{\scriptstyle\pi} \\
 Y & \xrightarrow{\fie} & X & = M/G
 \end{matrix}
 \label{eq:YtoM/G}
 \end{equation}
Here $M$ is a quasiprojective algebraic manifold with $K_M=0$ and $G$ a
finite auto\-morphism group of $M$ that acts trivially on a global basis
$s_M\in H^0(K_M)$. The object of study is the quotient variety $X=M/G$ and
its resolutions $Y\to X$, sometimes assumed to have $K_Y=0$. An important
motivating case is a finite subgroup $G\subset\SL(3,\C)$ acting on
$M=\C^3$.

 \subsection{Definition--Reassurance}\label{def:sX}
 The quotient varieties $X=M/G$ occuring here are singular. The theory of
minimal models of higher dimensional algebraic varieties (Mori theory) has
a whole battery of definitions that deal systematically with singular
varieties; here I only need one small item: the orbifolds $X$ here have
trivial canonical class $K_X=0$ (or trivial Serre--Grothendieck dualising
sheaf $\om_X=\Oh_X$). In concrete terms, this means the following: $X$ is
a complex $n$-fold (algebraic or analytic variety), nonsingular in
codimension~1, and its nonsingular locus $\NonSing X$ has an everywhere
nondegenerate holomorphic $n$-form $s_X$ (deduced from $s_M$). So $s_X$ is
a complex volume element at every nonsingular point of $X$, or in other
words, it is a global basis of $\Om^n_{\NonSing X}$. A resolution of
singularities $\fie\colon Y\to X$ is {\em crepant} if $K_Y=\fie^*K_X$ or
$\om_Y=\fie^*\om_X$, which simply means that $Y$ is a nonsingular $n$-fold
with $K_Y=0$ or $\om_Y=\Om^n_Y=\Oh_Y\cdot s_Y$, where $s_Y=\fie^*s_X$. More
generally, an arbitrary proper birational map $\fie\colon V\to X$ has a
{\em discrepancy divisor} $\De_\fie=\sum a_iE_i$ defined by
$K_V=\fie^*K_X+\sum a_iE_i$ with $a_i\ge0$; a divisor $E_i$ is {\em
crepant} if $a_i=0$. The discrepancy
$\De_\fie$ is the divisor of zeros on $V$ of the basic $n$-form $s_X$ on
$X$, generalising the divisor of zeros of the Jacobian determinant; in
Mori theory, it measures how far $V$ is from minimal.

 \subsection{Summary and slogan}
 I start with a preview of different approaches to the McKay
correspondence, which are treated in more detail in later sections. Each
of these approaches gives a result in the case of a finite subgroup
$G\subset\SL(3,\C)$ acting on $M=\C^3$.
 \begin{enumerate}
 \raggedbottom
 \item Gonzalez-Sprinberg and Verdier sheaves: the first direct link from
the representation theory of $G$ to the geometry of the resolution $Y\to X$
was the work of Gonzalez-Sprinberg and Verdier \cite{GSpV}: for a Kleinian
subgroup $G\subset\SL(2,\C)$, they constructed {\em sheaves $\sF_\rho$ on
$Y$, indexed by the irreducible representations of\/ $G$, whose first
Chern classes base the cohomology of\/ $Y$}.
 \item\label{it:str} String theory: the first hint of a McKay
correspondence in higher dimensions comes from work of the string
theorists Dixon, Harvey, Vafa and Witten \cite{DHVW} around 1985: if
$G\subset\SL(3,\C)$ and $Y\to X=\C^3/G$ is a crepant resolution of the
quotient $\C^3/G$, {\em the Euler number of\/ $Y$ equals the number of
conjugacy classes of\/ $G$ (or the number of its irreducible
representations).}
 \item Explicit methods: the finite subgroups $G\subset\SL(3,\C)$ are
classified, and work in the early 1990s of Roan, Ito, Markushevich and
others proved case-by-case {\em the existence of crepant resolutions, and
the validity of the formula of\/} \cite{DHVW} for the Betti numbers of $Y$.
 \item Valuation theory: for a finite subgroup $G\subset\SL(n,\C)$, the
paper \cite{IR} shows that $G$ has a grading by {\em age}, analogous to the
weight grading in Hodge theory, and proves that {\em the conjugacy classes
of\/ {\em junior} elements $g\in G$ (elements of age\/~$1$) correspond
one-to-one with the crepant divisors of a resolution} (more precisely,
their discrete valuations). This result holds for any $G\subset\SL(n,\C)$
and is intrinsic, classification-free; but for $n\ge4$ it only gives a
small part of a McKay correspondence (so far).
 \item Nakamura's $G$-Hilbert scheme: a resolution of singularities $Y\to
X$, even if it is a Mori minimal model theory, is not at all unique.
Moreover, if $X=M/G$, the construction of a resolution $Y$ need not have
much to do with $G$. In 1995, Nakamura made the revolutionary suggestion
that in many interesting cases, {\em the $G$-Hilbert scheme is a preferred
resolution\/ $Y$ of\/ $X$} (see \cite{IN2}, \cite{N}, \cite{R}). When this
holds, $Y$ is a ``very good'' moduli space over $M$, and the general yoga
of moduli suggests that there should be a ``tautological'' treatment of the
geometry of $Y$ (comparable to the cohomology of Grassmann varieties).
 \item\label{it:FM} Fourier--Mukai transform: the derived category $\D(V)$
of coherent sheaves on a variety $V$ (considered up to isomorphism of
triangulated categories) can be used as a geometric characteristic of $V$,
in place of K~theory or cohomology. The Fourier--Mukai transform is a
general method for constructing isomorphisms of derived categories (see
\cite{Mu}, \cite{O}, \cite{BO1}, \cite{Br}, \cite{BrM}). Bridgeland and
others \cite{BKR} have recently used this technique to prove that, if
$Y=\GHilb M$ is a crepant resolution, then $\D^G(M)=\D(Y)$. This implies
the corresponding result in K~theory.
 \item\label{it:moti} Motivic integration: the motivic integration of
Batyrev, Denef and Loeser, and Kontsevich is a rigorous and comparatively
simple mathematical trick that mimics some aspects of the path integrals
of QFT. Very roughly, if $\fie\colon Y\to X$ is a resolution of
singularities, possibly far from minimal, with discrepancy divisor
$K_Y-\fie^*K_X=\sum a_iE_i$, the calculation amounts to defining the {\em
stringy homology} of $X$ by picking only $\frac{1}{a_i+1}$th of the
homology of $E_i$. Quite remarkably, {\em this is well defined, agrees
with the predictions of \cite{DHVW} mentioned in (2) above, and provides
an exact form of the homological McKay correspondence for finite subgroups
$G\subset\SL(n,\C)$.}
 \item Explicit methods (bis): for a finite group $G\subset\SL(3,\C)$, the
results of (\ref{it:FM}) (maybe also (\ref{it:moti})) imply that
Gonzalez-Sprinberg--Verdier sheaves $\sF_\rho$ base the K~theory of the
resolution $Y\to X$, so that their Chern classes or Chern characters base
the cohomology. Reworking this in explicit terms presents a treasure chest
of delightful computational problems -- already the Abelian cases lead to
lovely pictures (compare \cite{R}, \cite{CR}, \cite{C2}).
 \end{enumerate}

I believe that many other approaches to the McKay correspondence remain to
be discovered, and many interrelations between the different approaches;
this problem area is recommended to afficionados of noncommutative
geometry, perverse sheaves, Gromov--Witten invariants, elliptic
cohomology, Chow groups, etc. Here is an attempt to describe the subject
in a single statement:

 \begin{princ}\label{princ:1.1} Let\/ $M$ be an algebraic manifold, $G$ a
group of automorphisms of\/ $M$, and\/ $Y\to X$ a resolution of
singularities of\/ $X$. Then the answer to any well posed question about
the geometry of\/ $Y$ is the\/ $G$-equivariant geometry of\/ $M$.
 \end{princ}

I give two illustrations
 \begin{enumerate}
 \item[I.] If $G\subset\SL(n,\C)$ acts on $\C^n$ and the quotient
$X=\C^n/G$ has a crepant resolution $Y\to X$, the homology or K~theory of
$Y$ is expected (or known) to be independent of $Y$. In this case, the
principle says that the homology or K~theory of $Y$ is the representation
theory of $G$ (equal to the $G$-equivariant geometry of $\C^n$ because
$\C^n$ is contractible).
 \item[II.] Let $M$ be a Calabi--Yau $n$-fold and $G$ a group of
automorphisms of $M$ that acts trivially on $\Om^n_M$. The stringy
homology of $X=M/G$ (see Sections~\ref{sec:dhvw} and~\ref{sec:mot_int})
is well defined by \cite{DL1}. The principle says that it must agree with
the $G$-equivariant homology of $M$. (I expand on what this means in
Section~\ref{sec:mot_int}.)
 \end{enumerate}

Viewed as an orbifold or stack, $X=M/G$ contains $M$ and the $G$ action,
and you can of course derive tautological question-and-answer pairs from
this (this is often popular as a source of questions after the talk). The
content of my slogan is that {\em the equivariant geometry of\/ $M$ already
knows about the crepant resolution\/ $Y\to X$}. Minimal models exist for
surfaces by classical work, and for 3-folds by Mori theory (or by explicit
methods). Minimal models of orbifolds by finite subgroups
$G\subset\SL(3,\C)$ provide infinitely many examples of local models of
Calabi--Yau 3-folds; calculating their Betti numbers or K~theory in a
priori terms is in no sense a tautology. If you prefer to think of the
singular $X$ as the fundamental object, and not resolve it (a perfectly
sensible alternative), the content is that {\em $X$ has invariants that
can be described from the orbifold $M/G$, but are birationally
invariant\/} under appropriate conventions about resolutions.

 \section{Age and discrepancy}\label{sec:age}
 Let $G\subset\SL(n,\C)$ be a finite group; any element $g\in G$ has
finite order $r$, say. For any such $r$, I choose at the outset a primitive
$r$th root of 1, say $\exp\frac{2\pi i}{r}$. A choice of eigenbasis
diagonalises the action of $g\in G$ on $M=\C^n$, giving
 \begin{equation}
 g=\diag(\ep^{a_1},\dots,\ep^{a_n})\quad\text{with $0\le a_i<r$.}
 \label{eq:ai}
 \end{equation}
I write $g=\frac{1}{r}(a_1,\dots,a_n)$, possibly depending on the choices
made. Now $\sum a_i\equiv 0$ mod~$r$ because $g\in\SL(n,\C)$. Following
\cite{IR}, define the {\em age} of $g$ by $\age g=\frac{1}{r}\sum a_i$. As
we will see, this is an analog of {\em weight} in Hodge theory; Denef and
Loeser \cite{DL2} refer to it by the long-winded but not inappropriate
term {\em valuation theoretic weight\/}. Clearly, $\age g$ is an integer
in the range $[0,\dots,n-1]$, and only the identity has age~0. The
elements of age~1 are {\em junior}.

Junior elements of $G$ give rise to crepant divisors of a resolution $V\to
M/G$ by the following toric mechanism (for more details, and a picture,
see \cite{IR}, 2.6--7). For $g\in G$ (not the identity), consider the
$a_i$ of (\ref{eq:ai}), and suppose $(a_1,\dots,a_n)\in\Z^n$ is primitive.
The coordinate subspace corresponding to the $x_i$ with $a_i=0$ is the
fixed locus $\Fix g$; it splits off as a direct product, and I assume that
all $a_i>0$ to short-cut some simple arguments. A useful example to bear
in mind is when all the $a_i=1$ (compare Example~\ref{exa:Vb}).

I view the integers $(a_1,\dots,a_n)$ as {\em weights}. They define the
grading $\wt x_i=a_i$ on the coordinate ring $\C[x_1,\dots,x_n]$, or
equivalently, the action $x_i\mapsto\la^{a_i}x_i$ of $\C^*$ on $M=\C^n$
that defines the weighted projective space
 \begin{equation}
 \PP(a_1,\dots,a_n)=(\C^n\setminus0)/\C^*.
 \notag
 \end{equation}
We obtain the weighted blowup $B_g\to M$ as the closed graph of the
quotient map $M\broken\PP(a_1,\dots,a_n)$; it has the exceptional divisor
$B_g\supset E_g=\PP(a_1,\dots,a_n)$. Obviously $g$ acts on $B_g$, and
fixes $E_g$ pointwise (because $g$ acts on $M$ as $\ep\in\C^*$). Therefore
$B_g\to B_g/\Span{g}$ is totally ramified along $E_g$.

 \begin{theo}[\cite{IR}, 2.6--7]
 Suppose that\/ $V\to X$ is any resolution of singularities of the
quotient\/ $X=M/G$. Then\/ $V$ contains a divisor\/ $F_g$ rationally
dominated by\/ $E_g$ under the rational map\/ $B_g\to M\broken V$. The
discrepancy of\/ $F_g$ is given by\/ $a_{F_g}=\age g-1$, and in particular
 \begin{equation}
 \text{$F_g$ is crepant} \iff \text{$g$ is junior.}
 \notag
 \end{equation}
Every crepant divisor of any resolution\/ $V$ occurs thus.
 \end{theo}

 \begin{pf}{Discussion of proof} Write $X_g=M/\Span{g}$ for the partial
quotient. Then $B_g/\Span{g}\to X_g$ is a partial resolution, with the
single exceptional divisor $E_g$. An easy toric calculation gives the
discrepancy of $E_g\subset B_g$ or $E_g\subset B_g/\Span{g}$ (compare
\cite{YPG}, 4.8): the standard basis of $\Om^n_M$ is $s_M=\dd
x_1\wedge\cdots \wedge\dd x_n$. For $K_{B_g}$, choose a Laurent monomial
$y_1=x^{\bm}$ of weight 1 (recall that the $a_i$ were coprime). Then $y_1$
is the defining equation of $E_g\subset B_g$ at a general point of $E_g$
(away from all the coordinate hyperplanes), and $y_1^r$ that of $E_g\subset
B_g/\Span g$. Choosing Laurent monomials $y_2,\dots,y_n$ forming a basis
of the lattice of monomials of weight 0, we get that
 \begin{equation}
 s_{B_g}=\dd y_1\wedge\cdots\wedge\dd y_n\in\Om^n_{B_g}
 \notag
 \end{equation}
 is the required basis. The discrepancy is the exponent $a$ in
$s_M=\text{(unit)}\cdot y_1^a s_{B_g}$, and is determined by weighty
considerations: $s_M$ has weight $\sum a_i$ and $s_{B_g}$ weight 1, so
$a=\sum a_i-1$. On the quotient $B_g/\Span g$ we only have $y_1^r$, so
we get the stated discrepancy $\frac{1}{r}\bigl(\sum
a_i-1-(r-1)\bigr)=\age g-1$.

The quotient morphism $M\to X$ is a Galois cover with group $G$; a cyclic
subgroup $\Span g$ corresponds to an intermediate cover $M\to
M/\Span{g}=X_g\to X$. The reduction to a cyclic group is in terms of
ramification theory; see \cite{IR}, 2.6--7. Roughly, over the general
point of any exceptional divisor $F$ of $V\to X$, the Galois extension of
function fields $k(X)\subset k(M)$ forms a tower, starting with a cyclic
ramified cover. For a crepant exceptional divisor, the cyclic ramification
can be chased back up to a conjugacy class of junior elements $g\in G$.
 \end{pf}

 \begin{rema} This argument works in all dimensions, but it only identifies
the {\em divisors} of a crepant resolution $Y$, and thus only gives a basis
of $H^2(Y,\Q)$ or $H_{2n-2}(Y,\Z)$ corresponding in McKay style to junior
conjugacy classes of $G$. In 3 dimensions, we used Poincar\'e duality to
bootstrap ourselves up to a basis of $H^*(Y,\Q)$ in \cite{IR}.
Historically, this was the first intrinsic proof of the conjectured
formula of \cite{DHVW} for the Betti numbers of a crepant resolution.

 As Brylinski \cite{B} remarks (following Mumford), if $V\to X$ is any
resolution, the group $G$ can be viewed as the fundamental group of $V$
minus the branch locus, so that an exceptional divisor $F$ of a resolution
$V$ corresponds directly to a conjugacy class of $G$ as a little
anticlockwise loop around $F$; for crepant divisors, this is of course the
same relation as in \cite{IR}. But I don't know how to use this idea to
get a well defined relation between, say, codimension 2 cycles of $Y$ and
age~2 conjugacy classes of $G$.
 \end{rema}

 \section{L'innomable}\label{sec:dhvw}
 This section is mainly for sociological and historical interest, but some
harmless hilarity may derive from my garrulous display of incompetence and
ignorance in physics.

 A theoretical prediction of string theory: Fermionic strings propagate in
10-dimensional space-time. Indeed, a universe of any other dimension would
have particles moving faster than the speed of light. Since this
prediction, on the face of it, contradicts the empirically observed
4-dimensions of space-time, string theorists want 6 of the dimensions to
be filled up with tiny Calabi--Yau 3-folds. (This means (i) a
6-dimensional Riemannian manifold with $\SU(3)$ holonomy, or (ii) a
complex manifold $V$ with a Ricci flat K\"ahler metric and $H^1(V,\R)=0$,
or (iii) an algebraic manifold $V$ with $K_V=0$ and $H^1(V,\Oh_V)=0$. It
seems that the holonomy or K\"ahler conditions on $V$, together with some
finite volume, are required by the physics, whereas making $V$
nonsingular, compact, and a constant fibre over macroscopic space-time are
just convenient choices when you try to guess a model.)

The two papers \cite{DHVW} were concerned with trying to calculate string
theory on examples of Calabi--Yau varieties obtained by dividing a
3-dimensional complex torus $M$ by a finite group $G$ preserving a basic
holomorphic 3-form, so that the stabiliser subgroup at any point is a
subgroup of $\SL(3,\C)$. A closed string on the quotient may lift either
to a closed string on the cover, or to a path that goes from $x$ to
$g\cdot x$. The latter are called {\em twisted sectors}. The physicists
need to take care of these in order to relate $\int_X$ to $G$-equivariant
$\int_M$, and they are the key to the form of the McKay correspondence in
Theorem~\ref{th:moti}, (4).

Taking limits is a tradition in physics, where the old is frequently the
limit of the new: Newtonian mechanics is the limit of special relativity
as $c\to\infty$, classical mechanics the limit of quantum physics as
$\hbar\to0$, groups and their Hopf algebras the limit of quantum groups as
$q\to1$. In string theory, if the scale (or radius of curvature) of the
tiny Calabi--Yau tends to zero, the theory should approximate ordinary
Lorentz 4-dimensional space-time, whereas letting it tend to macroscopic
proportions would approximate flat Lorentz 10-dimensional space-time. In
this context, the twisted sector near a point $x\in M^H$ plays the role of
strings that are topologically nontrivial, but are allowed to remain of
finite length (and so contribute to path integrals) as the scale becomes
large. To calculate something called the {\em $1$-loop partition
function}, DHVW considered mapping the elliptic curve $S^1\times S^1$
(with parameters $\si$ and $\tau$ along the copies of $S^1$) into $X$, or
the $\si,\tau$ square into $M$ with equivariant boundary conditions
depending on $g,h$. Thinking about twisted sectors and limits led DHVW (I
confess that their logic eludes me somewhat) to the formula
 \begin{equation}
 e\stringy(X)=e(M,G) :=\frac{1}{|G|}
 \sum_{\begin{smallmatrix} g,h\in G\\\text{commuting}\end{smallmatrix}}
 e(M^{\Span{g,h}}).
 \label{eq:DHVW}
 \end{equation}
Here $e(M,G)$ on the left-hand side is the {\em $G$-equivariant Euler
number of $M$}; on the right-hand side, the sum runs over all commuting
pairs of elements of $G$, $\Span{g,h}$ is the Abelian group they generate,
$M^{\Span{g,h}}$ its fixed locus in $M$, and $e$ is the usual Euler
number. The formula is a replacement for the Euler number of the singular
orbifold $X$. The papers \cite{DHVW} contain more-or-less explicitly the
conjecture that this number is the Euler number of a minimal resolution of
singularities.

It is not hard (see \cite{HH}, \cite{Roan} and \cite{Hwk}) to rearrange
the sums in (\ref{eq:DHVW}) to give
 \begin{equation}
 e\stringy(X)=e(M,G)
 = \sum_{[H]\subset G} e(X^H) \times \card\{[h]\in H\},
 \label{eq:HHRoan}
 \end{equation}
where (i) the first sum runs over conjugacy classes of subgroups $H\subset
G$; (ii) the stratum $X^H$ is the set of $x\in X$ such that $\Stab y$ is
conjugate to $H$ for any point $y\in M$ over $x$; (iii) the second factor
is the number of conjugacy classes in $H$. This means that $X^H\subset X$
contributes to $e(M,G)$ with multiplicity the representation theory of $H$.

 \begin{rema}
 The physicists want to do path integrals, that is, they want to integrate
some ``Action Man functional'' over the space of all paths or loops
$\ga\colon[0,1]\to Y$. This impossibly large integral is one of the major
schisms between math and fizz. The physicists learn a number of
computations in finite terms that approximate their path integrals, and
when sufficiently skilled and imaginative, can use these to derive
marvellous consequences; whereas the mathematicians give up on making
sense of the space of paths, and not infrequently derive satisfaction or a
misplaced sense of superiority from pointing out that the physicists'
calculations can equally well be used (or abused!)\ to prove $0=1$. Maybe
it's time some of us also evolved some skill and imagination. The motivic
integration treated in the next section builds a miniature model of the
physicists' path integral, by restricting first to germs of holomorphic
paths $\ga\colon U\to Y$, where $0\in U\subset\C$ is a neighbourhood of 0,
then to formal power series $\ga\colon\Spec\C\dbrk{z}\to Y$.
 \end{rema}

 \section{Motivic integration}\label{sec:mot_int}
 The material in this section is due to Batyrev \cite{Ba1}, \cite{Ba2},
Denef and Loeser \cite{DL1}, \cite{DL2} and Kontsevich \cite{K}. I
recommend Craw \cite{C1} as a readable first introduction to these ideas.

Rather than trying to restrict to crepant resolutions, take an arbitrary
normal crossing resolution $\fie\colon Y\to X$, marked by the discrepancy
divisor $D=\De_\fie=\sum_{i\in I} a_iD_i$ (here $I$ is the indexing set of
the components $D_i$). The normal crossing divisor $D$ defines a
stratification of $Y$, with
 \begin{equation}
 \text{{\em closed strata} $D_J=\bigcap_{j\in J} D_j$,}
 \quad \text{and} \quad \text{{\em open strata} $D\open_J=D_J\setminus
 \bigcup_{J'\supsetneq J} D_{J'}$}
 \notag
 \end{equation}
for $J\subset I$ (including, of course, $Y=D_\emptyset$ and $Y\setminus
D=D\open_\emptyset$).

Motivic integration is discussed and defined below, but it is convenient
to start from the answer: the {\em stringy motive} of $(Y,D)$, or of $X$
itself, turns out to be
 \begin{equation}
 h\stringy(X)=h(Y,D)=\sum_{J\subset I} [D\open_J]\cdot
 \prod_{j\in J} \frac{\LL-1}{\LL^{a_j+1}-1}\,.
 \label{eq:hstring}
 \end{equation}
Here $\LL=[\aff^1_\C]=[\C]$ is the Tate motive, and the formula takes place
in a certain ring of motives with formal power series in $\LL\1$ adjoined.
We will worry about the coefficient ring later, but in lucky cases it will
happen that the cyclotomic polynomials in the denominators cancel out,
leaving an integral motive (see Example~\ref{exa:Vb} and \cite{Hwk} for
examples). It follows from Theorem~\ref{th:moti}, (2) and~(3) that
$h(Y,D)$ is independent of the choice of the normal crossing resolution
$Y$, so depends only on $X$. In the case when $D=aE$ has a single
component with discrepancy $a$, it boils down to
 \begin{equation}
 [Y-E]+\frac{[E]}{1+\LL+\LL^2+\cdots+\LL^a} =
 [Y-E]+\frac{[E]}{[\PP^a]}.
 \label{eq:1/Pa}
 \end{equation}

 \begin{exem}\label{exa:Vb}
Let $n=ab$, and consider the $n$-fold quotient singularity $X$ of type
$\frac{1}{b}(1,\dots,1)$, that is, the quotient $\C^n/(\Z/b)$, with the
diagonal action of $\ep=\exp\frac{2\pi i}{b}$. It is the cone over the
$b$th Veronese embedding of $\PP^{n-1}$, so that its resolution $Y\to X$
has exceptional divisor $E=\PP^{n-1}$ with $\Oh_E(E)=\Oh_{\PP^{n-1}}(-b)$.
The discrepancy is $a-1$, to fit the adjunction formula, with
$K_Y=(a-1)E$, and $K_E=\Oh_E(aE)=\Oh_{\PP^{n-1}}(-n)$.

Now whereas $Y$ is homotopy equivalent to $\PP^{n-1}$, so has $n$
homology classes, one in each dimension $0,2,\dots,2(n-1)$, the effect
of dividing by $[\PP^{a-1}]$ in (\ref{eq:1/Pa}) is to throw away most of
these, leaving only the $b$ stringy homology classes in dimension
$0,2a,4a,2a(b-1)$. This is exactly what we need for the McKay
correspondence: the $b$ elements of $\Z/b$ have age $0,a,2a,\dots,a(b-1)$
and correspond to the stringy classes in dimension $2ia$.
\end{exem}

\begin{exem}
Consider the blowup $\si\colon Y_1\to Y$ of a subvariety $C\subset Y$ that 
intersects all the strata of $D$ transversally, and set
$D_1=\si^*D+(c-1)E$, where $E=\si\1C$ is the exceptional divisor of the
blowup and $c=\codim C$. The coefficient is the discrepancy of $E$, so
that $K_{Y_1}-D_1=\si^*(K_Y-D)$. It is an exercise to see that
 \[
 h(Y,D)=h(Y_1,D_1).
 \]
(This is rather trivial if $C\cap D=\emptyset$ in view of Grothendieck's
{\em formule clef\/} for the motive of a blowup; see \cite{Hwk} for more
hints.) This is good evidence for the birational invariance of $h(Y,D)$.
\end{exem}

I now describe briefly the mechanics of motivic integration, following
\cite{C1}. Start from the Grothendieck ring $K_0(\sV)$ of classes of
varieties under the equivalence relation $[V]=[V\setminus W]+[W]$.
Addition and multiplication are quite harmless. The Tate motive is
$\LL=[\aff^1_\C]=[\C]$. We formally adjoin $\LL\1$ to $K_0(\sV)$, and make
a fairly mild $(\LL\1)$-adic completion to give the value ring
$R=\widehat K_0(\sV)[\LL\1]$. This value ring is the really clever thing
about the whole construction. (Exercise: $(\LL^a-1)\1$ can be written as a
formal power series in $\LL\1$, so all the terms on the right-hand side of
(\ref{eq:hstring}) are in $R$.)

Motivic integration takes place over the infinite jet space $J_\infty Y$,
which coincides with the set $Y(\C\dbrk{z})$ of points of $Y$ with values
in the formal power series ring $\C\dbrk{z}$. An element $\ga\in
Y(\C\dbrk{z})$ is a point $y=\ga(0)\in Y$ together with a formal arc
$\ga\colon\Spec\C\dbrk{z}\to Y$ starting at $y$; if convergent, $\ga$ is
the Taylor series of a holomorphic germ $\wave\ga\colon(\C,0)\to Y$. The
infinite jet space $J_\infty Y$ is the profinite limit $\varprojlim_k
J_kY$ of the finite jet spaces $J_k Y$; recall that $J_0Y=Y$, $J_1Y$ is
the total space of the tangent bundle $T_Y$, and $J_{k+1}\to J_k$ is a
$\C^n$-fibre bundle.

The projection maps $\pi_k\colon J_\infty Y\to J_k$ of the profinite limit
allow us to define a {\em cylinder set\/} in $J_\infty Y$ to be
$\pi_k\1(B_k)$ for a constructible set $B_k\subset J_k$. The measure on
$J_\infty Y$ is initially defined on these, by setting\footnote{The papers
\cite{DL1} and \cite{C1} have the exponent $\LL^{-n(k+1)}$. This is just a
normalising convention, giving $h(Y,D)=[Y]\cdot\LL^{-n}$ in
Theorem~\ref{th:moti}, (1), and making the motive of $Y$ 0-dimensional. I
prefer my version.\label{ft:Ln}}
 \begin{equation}
 \mu(\pi_k\1(B_k)):=[B_k]\cdot\LL^{-nk}\in R.
 \label{eq:mu}
 \end{equation}
It is straightforward to see that this is independent of $k$, and is a
``finitely additive measure''.

As our measurable functions, consider an effective divisor $D$ on $Y$, and
define a function $F_D\colon J_\infty Y\to\Z_\ge0$ by $F_D(\ga)=D\cdot\ga$
(intersection number). In other words, suppose $\ga(0)=P\in Y$ and let
$g_D$ be the local defining equation of $D$ at $P$; then $F_D(\ga)$ is the
order of $\ga^*(g_D)\in\C\dbrk{z}$. Since the first $s$ coefficients of
$\ga^*(g_D)$ clearly only depend on $\pi_s(\ga)\in J_s$, it is obvious that
$F_D\1(s)$ is a cylinder set.

The grand definition is now: for $Y$ a nonsingular variety and $D$ a
normal crossing divisor, the motivic integral is
 \begin{equation}
 h(Y,D)=\int_{J_\infty Y}\LL^{-F_D}:=
 \sum_{s\in\Z_{\ge0}}\mu\bigl(F_D\1(s)\bigr)\cdot\LL^{-s}\in R.
 \label{eq:hYD}
 \end{equation}

 \begin{rema}
 I omit some tricky details on convergence required to get a genuine
measure (involving the $[\LL\1]$-adic completion). To tell the truth, I
don't know if they are at all essential. A basic point for applications is
that the measure of $F_D\1(s)$ tends to~0 as $s\to\infty$; this is
plausible enough (because arcs $\ga$ with $\ga\cdot D\ge s$ have
codimension $\ge s$ in $J_\infty Y$), and is an intuitive reason behind
birational invariance: the arcs in a Zariski closed subset of $Y$ have
measure zero.
 \end{rema}

 \begin{theo}\label{th:moti} $h(Y,D)$ of (\ref{eq:hYD}) has the following
properties:
 \begin{enumerate}
 \item If $D=0$ then $h(Y,D)=[Y]$.
 \item $h(Y,D)$ is calculated by the right-hand side of (\ref{eq:hstring}).
 \item Birational invariance: let $Y',D'$ and $Y,D$ be pairs, and
$\fie\colon Y'\to Y$ a birational morphism such that
$K_{Y'}-D'=\fie^*(K_Y-D)$; then
 \begin{equation}
 h(Y',D')=h(Y,D).
 \notag
 \end{equation}
 \item If\/ $X=M/G$ is as in Assumption~\ref{ssec:ass}, $Y\to X$ a normal
crossing resolution, and $D$ the discrepancy, then
 \begin{equation}
 h\stringy(X)=h(Y,D)
 =\sum_{[H]\subset G} \left[X^H\right] \cdot\sum_{[g]\in H}\LL^{\age g},
 \label{eq:McK}
 \end{equation}
where the range of summation is as in (\ref{eq:HHRoan}), and the second
sum is over conjugacy classes in $H$.
 \end{enumerate}
 \end{theo}

 \begin{pf}{Discussion of proof}
 I give some indications, leaving most of the proof as references to
\cite{DL1} and~\cite{DL2}. Alternatively, do them as exercises (see
\cite{Hwk} for more hints). The key point of the proof is that, whatever
its substance, (\ref{eq:hYD}) has the formal properties of an integral,
and is subject to the same kind of {\em change of variables\/} formula. In
the words of the Master:
 \begin{quote}
``La th\'eorie consiste pour l'essentiel dans des questions de variance''
 \end{quote}
(\cite{H}, Introduction). Note first that the condition in (3) says that
$D'-D=\div(\Jac\fie)$ is the divisor of zeros of the Jacobian determinant
of $\fie$ (I omit $\fie^*$ from now on). Composition defines a map
$j_\fie\colon J_\infty Y'\to J_\infty Y$, and, unless it falls entirely in
the locus of indeterminacy of $\fie\1$, an arc in $Y$ has a birational
transform as an arc in $Y'$; in other words, away from subsets of measure
zero, $j_\fie$ is a bijection on the infinite jet spaces. For (3), it
remains only to stratify the finite jet spaces $J_kY'$ and $J_kY$ so that
the corresponding morphism $j_k\colon J_kY'\to J_kY$ is a $\C^t$-bundle on
each stratum with $F_{D'-D}(\ga)=\div(\Jac\fie)\cdot\ga=t$ (see \cite{DL1},
Lemma~3.4 and \cite{Hwk}).

 (2) is proved in \cite{DL1}, Proposition~6.3.2, \cite{Ba2}, Theorem~6.28,
and worked out in detail in \cite{C1}, Theorem~1.16. The proof of (4)
consists of two steps, relating to the two morphisms $\pi\colon M\to X$ and
$\fie\colon Y\to X$ of Assumption~\ref{ssec:ass}.

\begin{pf}{Step~I} We translate the twisted sectors of \cite{DHVW} into
the language of formal arcs, obtaining the stratification (\ref{eq:tw})
below.

Let $y\in M^H$ be a point with $\Stab y=H$ and $x=\pi(y)\in X^H$. As at
the start of Section~\ref{sec:age}, suppose that $r$ is an integer
divisible by the order of each $g\in H$, and choose an $r$th root $\ep$ of
1 and an $r$th root $\ze=z^{1/r}$ of the parameter used for formal arc, so
that a formal arc $\ga$ at $x\in X$ parametrised by $z$ lifts to a formal
arc at $y\in M$ parametrised by $\ze$. Unless $\ga$ falls entirely in the
branch locus of $\pi\colon M\to X$, there is a unique conjugacy class
$g\in H$ defined by $\ga(\ep\ze)=g\ga(\ze)$. Here $g$ is the {\em twisted
sector}, the conjugacy class of $\ga$ in the local fundamental group $H$
(where $\ga$ is viewed as a little loop in $X$ minus the branch locus).

This argument shows that, after we delete the subset of arcs falling
entirely in the branch locus (which has infinite codimension, so measure
zero) the infinite jet space $J_\infty X$ is a disjoint union
 \begin{equation}
 J_\infty X
 =\coprod_{[H]\subset G} \quad \coprod_{[g]\in H} \quad J_\infty^{H,g} Y,
 \label{eq:tw}
 \end{equation}
where $H,g$ are as in (\ref{eq:HHRoan}), and $J_\infty^{H,g} Y$ is the set
of arcs with $\ga(0)\in X^H$ in the twisted sector $g$.
\end{pf}

\begin{pf}{Step~II} Using change of variables as in the proof of (3), one
calculates that $J_\infty^{H,g} Y$ contributes $X^H\cdot\LL^{g\1}$ to
$h(Y,D)$ (\cite{DL2}, Lemma~4.3). The difference in appearance of the
formulas here and in \cite{DL2} are explained by two trivial shifts of
notation: as explained in the footnote on page~\pageref{ft:Ln}, my measure
is $\LL^n$ times theirs; and they diagonalise $g$ as $\ep^{e_i}$ with
$1\le e_i\le r$, defining $w(g)=\frac{1}{r}\sum e_i=n-\age(g\1)$.
 \end{pf}
 \end{pf}

 \begin{rema}\label{rem:Dvol}
 Statement (4) is an exact analogue of the \cite{DHVW} formula
(\ref{eq:HHRoan}), saying that the stratum $X^H$ appears in the stringy
homology of $Y$ multiplied by the set of conjugacy classes in $H$.

As discussed in Definition~\ref{def:sX}, the discrepancy $D=\div s_X$ is
the divisor of zeros of $s_X$, the global basis of $\Om^n_{\NonSing X}$.
In the normal course of events, integrating {\em functions} on $Y$
requires a volume form; here we take $s_X$ as a holomorphic volume form,
viewing its zeros on $D$ as scaling down the contribution from
neighbourhood of the discrepant exceptional divisors. This is what
produces a birationally invariant answer.
 \end{rema}

 \section{Hilbert schemes of $G$-orbits}
 This section explains the definition of the $G$-orbit Hilbert scheme
$\GHilb M$, and Nakamura's idea of using it to resolve certain quotient
singularities. We know by general results (especially Hironaka's
theorems) that the singularities of a quotient variety $X=M/G$ can be
resolved somehow-or-other, but the construction of an actual resolution is
messy, involves lots of choices, and will probably have almost nothing to
do with the group action. Around 1995, Ito and Nakamura observed that in
the case of $G\subset\SL(2,\C)$, {\em the Hilbert scheme\/ $\GHilb\C^2$
of\/ $G$-clusters is a crepant resolution of the quotient\/ $\C^2/G$.}
Nakamura conjectured that this continues to hold for $G\subset\SL(3,\C)$,
and this has since been confirmed and extended to some other cases by work
of Bridgeland and others (see \cite{BKR} and Theorem~\ref{th:BKR}).

First, a {\em cluster} in a variety $M$ (say, quasiprojective and
nonsingular) is a 0-dimensional subscheme $Z\subset M$, defined by an ideal
$\sI_Z\subset\Oh_M$, so that the cokernel $\Oh_Z=\Oh_M/\sI_Z$ is a finite
dimensional $\C$-vector space. The {\em degree} of $Z$ is the dimension of
$\Oh_Z$. Like the intersection of two plane curves in Bezout's theorem, a
cluster $Z$ may consist of reduced points $Z=P_1+\cdots+P_N$, or may have a
nonreduced structure; in the latter case, we {\em keep track of the ideal\/
$\sI_Z\subset\Oh_M$}, as a way of using algebraic equations to keep
information about the relative positions when some of the points $P_i$
come together. For example,
 \[
 (x^2,xy,y^2) \quad\text{and}\quad (x-ay-by^2,y^3) \quad
 \text{for any $a,b\in\C$}
 \]
are clusters of degree 3 supported at $0\in\C^2$.

 \begin{lem}
 All clusters $Z\subset M$ of given degree $N$ in $M$ are parametrised by
a quasi\-projective scheme $\Hilb^NM$, which is a fine moduli space.
 \end{lem}

 \begin{proof} The assertion is quite elementary. $M$ is quasiprojective;
choose an embedding $M\subset\PP^s$. Every ideal $\sI_Z\subset\Oh_M$ of
codimension $N$ defines and is defined by a codimension $N$ vector
subspace
 \[
 H^0(\PP^s,\sI_Z(N))\subset H^0(\PP^s,\Oh_{\PP^s}(N)),
 \]
the forms of degree $N$ vanishing on $Z$ (same $N$). Subspaces of given
codimension are parametrised by a Grassmann variety, and the condition that
a space of forms defines a cluster of degree $N$ in $M$ is a locally closed
condition. (It can be written in terms of rank of a matrix $=N$.)
 \end{proof}

 \begin{rema}
 The map $\Hilb^NM\to S^NM$ to the symmetric product, defined at the level
of sets by $Z\mapsto\Supp Z$, is a morphism of schemes, the Hilbert--Chow
morphism (see \cite{GIT}, Chapter~5, \S4). For a curve, $\Hilb^NC$ is just
the symmetric product $S^NC$, which is itself already nonsingular. For a
surface, the symmetric product $S^NS$ is singular at the diagonals, and
$\Hilb^NS\to S^NS$ is a crepant resolution, in fact, a symplectic
resolution; see \cite{IN2}, \S6. But $\Hilb^NM$ is singular as soon as
$\dim M\ge3$ and $N=\deg Z\ge4$, and usually even has components of excess
dimension.
 \end{rema}

 \begin{prodf}[Ito and Nakamura]
 Let\/ $G$ be a finite group of order\/ $N$ acting faithfully on an
algebraic manifold\/ $M$; consider the action of\/ $G$ on\/ $\Hilb^NM$ and
its fixed locus\/ $(\Hilb^NM)^G$. This has a unique irreducible component
that contains a general orbit\/ $G\cdot y$ of\/ $G$ on\/ $M$. This
component is defined to be the\/ {\em $G$-Hilbert scheme}, and denoted
by\/ $\GHilb M$. The composite\/ $\GHilb M\into\Hilb^NM\to S^NM$ induces a
Hilbert--Chow morphism\/ $\GHilb M\to M/G$ which is proper and birational.

A cluster\/ $Z\in\GHilb$ is\/ $G$-invariant, and is called a\/ {\em
$G$-cluster}; its defining ideal\/ $\sI_Z$ is\/ $G$-invariant, and as a
representation of\/ $G$, the quotient\/ $\Oh_Z=\Oh_M/\sI_Z$ is the regular
representation\/ $\C[G]$.

\rm See also \cite{CR}, 4.1 for a rival definition and a comparison between
the two.

 \end{prodf}

 \begin{proof}
 The general orbit $G\cdot y$ consists of $N$ points permuted simply
transitively by $G$, so is a $G$-invariant cluster in $(\Hilb^NM)^G$.
These orbits fill out an irreducible open set in $(\Hilb^NM)^G$, because a
small $G$-invariant deformation of $G\cdot y$ is clearly still a set of
$N$ distinct points permuted by $G$ and disjoint from any fixed locus. The
closure of this component is $\GHilb M$ by definition. The composite
$\GHilb M\into\Hilb^NM\to S^NM$ is a morphism; by definition, a dense open
set of $\GHilb M$ consists of general orbits $G\cdot y$, and these maps to
orbits in $S^NM$, that is, to $M/G$.

Finally, the quotient sheaves $\Oh_Z$ for $Z\in\GHilb M$ fit together as a
locally free sheaf $\Oh_{\sZ}$ over $\GHilb M$, with a $G$-action that
makes it the regular representation on a dense open set. Its isotypical
decomposition under the idempotents of $\C[G]$ is a direct sum, so each
component must also vary as a locally free sheaf, therefore
$\Oh_Z\iso\C[G]$ for every $Z\in\GHilb M$ (since $\GHilb M$ is defined to
be irreducible).
 \end{proof}

The $G$-Hilbert scheme is a crepant resolution for finite groups
$G\subset\SL(3,\C)$. The general case of this is proved by Bridgeland and
others \cite{BKR} using derived category methods and a homological
characterisation of regularity. For a diagonal Abelian group, $\AHilb\C^3$
is a completely explicit construction of Nakamura (see \cite{N} and
\cite{CR}): the monomial $xyz$ is $A$-invariant, and every $G$-cluster $Z$
is defined by 7 (possibly redundant) equations of the form
 \[
 \begin{array}{c}
 x^{a+1}=\la y^dz^g\\
 y^{b+1}=\mu z^ex^h\\
 z^{c+1}=\nu x^fy^i
 \end{array}\quad
 \begin{array}{c}
 y^{d+1}z^{g+1}=\al x^a\\
 z^{e+1}x^{h+1}=\be y^b\\
 x^{f+1}y^{i+1}=\ga z^c
 \end{array}\quad\text{and}\quad
 xyz=\xi
 \]
for appropriate exponents $a,\dots,i$ and coefficients $\al,\dots,\xi$
satisfying $\al\la=\be\mu=\ga\nu=\xi$. The monomial basis of $\Oh_Z$ forms
a tripod shaped Newton polygon in the plane lattice $\Z^2$ of Laurent
monomials modulo $xyz$; this lattice is naturally the universal cover of
the McKay quiver and the tripod is a choice of fundamental domain for the
covering group (see \cite{N} and \cite{R} for pictures). The explicit
calculations remain an interesting challenge in the non-Abelian cases,
e.g., in the trihedral case.

 \begin{exem}\label{exa:4f}
 These results are known to fail for finite $G\subset\SL(4,\C)$. In the
first place, most quotient singularities $X=\C^4/G$ do not have any crepant
resolution. For example, the series of cyclic quotient singularities
$\C^4/(\Z/r)$ of type $\frac{1}{r}(1,r-1,i,r-i)$ have no junior elements,
so are terminal; compare Example~\ref{exa:Vb}. These examples motivated
the initial exploration of stringy homology in \cite{BD}.

 Next, even when a crepant resolution exists, the $G$-Hilbert scheme may be
singular or discrepant or both. A simple example is the quotient
singularity $\C^4/G$ by the maximal diagonal subgroup
$(\Z/2)^{\oplus3}\subset\SL(4,\C)$ of exponent 2. The junior simplex $\De$
has all the midpoints of the edges $\frac{1}{2}(1,1,0,0)$ etc., as lattice
points. This has several subdivisions into basic simplexes, giving crepant
resolutions, but none that is symmetric under permuting the coordinates --
the only symmetric thing you can do is chop off the 4 basic simplexes at
the corners, leaving a terminal simplex of volume 2. On the other hand,
$\GHilb\C^4$ is obviously symmetric.
 \end{exem}

 \section{Coherent derived category}\label{sec:der_cat}

Grothendieck and Verdier introduced the derived category $\D(X)$ of
coherent sheaves on a variety $X$ in the 1960s as a technical convenience
in homological algebra; it has enjoyed an unfortunate reputation for
technicality and abstraction ever since then. Recently, however, it has
been increasingly used as a geometric characteristic of $X$ similar to
K~theory: whereas K~theory works with the group of bundles or sheaves
modulo the relation $F=F'+F''$ for every short exact sequence $0\to F'\to
F\to F''\to0$, the derived category $\D(X)$ consists of complexes
$F\hidot$ modulo the relation of quasi-isomorphism (defined at the start
of the theory, and thankfully never referred to again). Following Mukai's
pioneering work \cite{Mu} for Abelian varieties, Orlov and Bondal
\cite{O}, \cite{BO1} have advocated the idea of considering the derived
category $\D(X)$ (up to isomorphism of triangulated categories) as a
geometric characteristic of $X$. {From} this point of view, $\D(X)$
behaves like an enriched version of K~theory.

A variety $X$ with $\pm K_X$ ample can be reconstructed from its derived
category $\D(X)$ (as a triangulated category) \cite{BO1}, but if $K_X=0$
(notably for an Abelian variety or a K3 surface), the same triangulated
category may occur as $\D(X)$ for different $X$, or there may be
infinitely many symmetries of $\D(X)$ not arising from automorphisms of
$X$. Isomorphisms $\D(X)\iso\D(Y)$ arise as {\em Fourier--Mukai transforms
$\Phi^{\sF}_{X\to Y}$} corresponding to a sheaf $\sF$ on $X\times Y$,
defined as the composite of the functors $p_X^*$, $\tensor\sF$ and
$q_{Y*}$ (more precisely, their derived functors); for an up-to-date
treatment, see \cite{Br} and the references given there. In practice, $Y$
is most frequently a moduli space of coherent sheaves on $X$ and $\sF$ the
universal sheaf over $X\times Y$, so that $Y$ parametrises sheaves $F_y$
on $X$; in very good cases, the apparatus of moduli functors, stable
bundles, and deformation theory gives essentially for free that the $F_y$
have orthonormality properties under $\Ext$ functors (formally analogous to
those of trig functions in the theory of Fourier transform).

Let $M$ be a nonsingular quasiprojective $n$-fold with $K_M=0$, and $G$ a
finite group acting on $M$, with trivial action on $K_M$. Set $Y=\GHilb
M$. Since $Y$ is a fine moduli space for $G$-clusters $Z\subset M$, there
is a universal $G$-cluster $\sZ\subset Y\times M$, fitting in a diagram
 \begin{equation}
 \renewcommand{\arraystretch}{1.2}
 \begin{matrix}
 \sZ & \xrightarrow{q} & M \\
 {\scriptstyle{p}} \big\downarrow \hphantom{\scriptstyle{p}}
 && \hphantom{\scriptstyle{\pi}}\big\downarrow{\scriptstyle{\pi}} \\
 Y & \xrightarrow{\fie} & X
 \end{matrix}
 \label{eq:FM}
 \end{equation}
Bridgeland and others \cite{BKR} prove the following theorem.

 \begin{theo}\label{th:BKR}
 Suppose that the inverse image of the diagonal\/
$(\fie\times\fie)\1(\De_X)$ has dimension\/ $\le n+1$ (automatic for
$n=3$). Then\/ $Y$ is a crepant resolution of\/ $X$ and the Fourier--Mukai
functor $\Phi=\R q_*\circ p^*\colon\D(Y)\to\D^G(M)$ is an equivalence of
categories.
 \end{theo}

Once we know that $Y$ is a crepant resolution, $\om_M$ is trivial as a
$G$-sheaf and $\om_Y$ is trivial, so that both the derived categories
$\D^G(M)$ and $\D(Y)$ have Serre duality functors; the remainder of the
proof is then standard Fourier--Mukai technology. However, the surprising
thing here is Bridgeland's derivation of the nonsingularity of $Y$ from
the famous theorem of commutative algebra known for a long time as Serre's
``Intersection conjecture''.

 \section{Fin de partie}

 Samuel Beckett's play of the same title has the wonderful line:
 \begin{quote}
 ``Personne au monde n'a jamais pens\'e aussi tordu que nous.''
 \end{quote}
This seems to reflect a truth about math research: progress beyond the
obvious takes really twisted thinking. In this spirit, let me raise all
the open questions I can think of.

There are two basic flavours of McKay correspondence:
 \begin{enumerate}
 \item conjugacy classes of $G$ $\bij$ homology of $Y$ (or stringy
homology); and
 \item representations of $G$ $\bij$ derived category $\D(Y)$ or K~theory
of
$Y$.
 \end{enumerate}
 Is there a ``bivariant'' version of the correspondence containing both
(1) and (2) at the same time? For example, in some contexts, $\sD$-modules
or perverse sheaves manage to accommodate both coherent and topological
cohomology. Note that (1) and (2) achieve a well posed question in
completely different ways: (1) takes accounts of discrepancy
systematically, whereas (2) currently only works under the very strict
condition that $Y=\GHilb$ is a crepant resolution.

The representation theory of finite groups has two ingredients, conjugacy
classes and irreducible representations, and a character table, which is a
nonsingular matrix making them ``dual'' (I apologise to group theorists
for this gratuitous vulgarity). Although in substance very different, the
homology and K~theory of a variety $Y$ could be described in similar
terms. In cases when McKay holds, is there any direct relation?

All the different approaches to McKay described here have one thing in
common: none of them seems to say anything very useful about multiplicative
structures. The following questions seem most likely to be approachable:
can tensor product of $G$-modules and tensor product in K~theory of $Y$ be
related? Can you reconstruct the McKay quiver in $\D(Y)$ or $K_0Y$?

Motivic integration takes a fraction of the homology of a discrepant
exceptional divisor, say, half the homology of the exceptional $\PP^3$ for
the quotient singularity $\C^4/(\Z/2)$ (the cone on the second Veronese
embedding $v_2(\PP^3)$). In contrast, half of a derived category is
something no-one has ever seen. In the case of $v_2(\PP^3)$, the
Gonzalez-Sprinberg--Verdier sheaves corresponding to the characters $\pm1$
are $\Oh_Y$ and $\Oh_Y(1)$.  Breaking up the derived category $\D(Y)$ into
two bits, one of which will correspond to the representations of $\Z/2$,
doesn't seem to make any sense. On the other hand, in this case we can
extend the action of $\Z/2$ to the action $\frac{1}{4}(1,1,1,1)$ of
$\Z/4$, whose quotient does have a crepant resolution.

Another general problem area: resolutions of Gorenstein quotient
singularities give a collection of examples of Calabi--Yau 3-folds with
very nice properties: the homology of the resolution is well defined
(independent of the choice of resolution), and the homology and K~theory
are closely related by something like a duality. Do these properties hold
for Calabi--Yau 3-folds more generally? It seems very likely that
birational Calabi--Yau 3-folds have isomorphic derived categories, but so
far this only seems to be established when they are related by classic
flops \cite{BO2}.

Part of motivic integration is the simple idea of using $\fie^*s_X$ as the
volume form, even though it vanishes along the discrepancy divisor $D$ 
(compare Remark~\ref{rem:Dvol}). Maybe this idea can be used with
differentials on $X$ itself (not passing to $J_\infty X$) to get
birationally invariant de Rham and Hodge cohomology?

Elliptic cohomology is another area of geometry with an alleged stringy
interpretation -- as the index of the Dirac operator on the space of
loops. Could part of this theory have a rigorous treatment in terms of
spaces of formal arcs, like motivic integration in
Section~\ref{sec:mot_int}? If we believe that the elliptic cohomology of
$M/G$ has a well defined answer (see Totaro \cite{T} for some evidence)
then Principle~\ref{princ:1.1} predicts what the answer must look like in
a whole pile of substantial cases.

Which Gorenstein quotient singularities admit crepant resolutions? Since
4-fold singularities usually do not have crepant resolutions, those that
do are of particular interest; see \cite{DHZ} for examples. How does this
relate to complex symplectic geometry? The papers of Verbitsky \cite{Vb}
and Kaledin \cite{Ka1}, \cite{Ka2} study crepant resolutions and related
issues for symplectic quotient singularities. When crepant resolutions
exist they are symplectic \cite{Vb}, therefore ``semismall'', giving a
complete and elegant solution to the homological form (1) of the McKay
correspondence \cite{Ka2}. Is it possible that there is a ``special''
geo\-metry in 3 complex dimensions (such as complexified imaginary
quaternions), like symplectic or hyper-K\"ahler geometry for complex
surfaces or 4-folds, that explain why crepant resolutions exist for
3-folds?

How should we interpret Nakamura's results and conjectures on $\GHilb$? If
a crepant resolution exists, it would be exceedingly convenient to be able
to describe it as a fine moduli space of something; $G$-clusters have no
especially privileged role, but the requirement that the space be
birational to $M/G$ seems to impose some relation with the moduli space of
group orbits. Nakamura and Nakajima have raised the question of whether
the {\em other} crepant resolutions (after a flop) can also be interpreted
as moduli, for example as Quot schemes; a single convincing example of
this would add weight to their suggestion. Do the crepant resolutions in
Example~\ref{exa:4f} have interpretations as moduli?

\end{document}